\documentclass{article}

\usepackage{
 amsmath,
 amsxtra,
 amsthm,
 amssymb,
 etex,
 mathrsfs,
 }
\usepackage[all]{xy}

\newtheorem{theorem}{Theorem}[section]
\newtheorem{lemma}[theorem]{Lemma}

\newtheorem{proposition}[theorem]{Proposition}
\newtheorem{corollary}[theorem]{Corollary}

\theoremstyle{definition}

\newtheorem{defn}[theorem]{Definition}
\newtheorem{remark}[theorem]{Remark}

\newcommand{\bd}{\begin{defn}}
\newcommand{\ed}{\end{defn}}
\newcommand{\bl}{\begin{lemma}}
\newcommand{\el}{\end{lemma}}
\newcommand{\bp}{\begin{proposition}}
\newcommand{\ep}{\end{proposition}}
\newcommand{\bt}{\begin{theorem}}
\newcommand{\et}{\end{theorem}}
\newcommand{\bc}{\begin{corollary}}
\newcommand{\ec}{\end{corollary}}
\newcommand{\br}{\begin{remark}}
\newcommand{\er}{\end{remark}}
\newcommand{\ba}{\begin{array}}
\newcommand{\ea}{\end{array}}
\newcommand{\bpf}{\begin{proof}}
\newcommand{\epf}{\end{proof}}

\newcommand{\Q}{\mathbb{Q}}

\newcommand{\Z}{\mathbb{Z}}
\newcommand{\Zp}{\mathbb{Z}_p}
\newcommand{\Op}{\mathcal{O}}

\newcommand{\e}{\varepsilon}
\newcommand{\Ga}{\Gamma}

\newcommand{\la}{\lambda}
\newcommand{\si}{\sigma}
\newcommand{\om}{\omega}
\newcommand{\ze}{\zeta}

\DeclareMathOperator{\Gal}{Gal}

\DeclareMathOperator{\Cl}{Cl}

\newcommand{\lra}{\longrightarrow}

\newcommand{\ot}{\otimes}

\newcommand{\ilim}{\displaystyle \mathop{\varinjlim}\limits}
\newcommand{\plim}{\displaystyle \mathop{\varprojlim}\limits}

\numberwithin{equation}{section}

\begin{document}
\title{On the structure of even $K$-groups of rings of algebraic integers}
 \author{
  Meng Fai Lim\footnote{School of Mathematics and Statistics $\&$ Hubei Key Laboratory of Mathematical Sciences,
Central China Normal University, Wuhan, 430079, P.R.China.
 E-mail: \texttt{limmf@ccnu.edu.cn}} }
\date{}
\maketitle

\begin{abstract} \footnotesize
\noindent In this paper, we describe the higher even $K$-groups of the ring of integers of a number field in terms of class groups of an appropriate extension of the number field in question. This is a natural extension of the previous collective works of Browkin, Keune and Kolster, where they considered the case of $K_2$.  We then revisit the Kummer's criterion of totally real fields as generalized by Greenberg and Kida. In particular, we give an algebraic $K$-theoretical formulation of the said criterion which we will prove using the algebraic $K$-theoretical results developed here. 


\end{abstract}

\section{Introduction}

Throughout, $p$ will always denote a fixed odd prime. Let $F$ be a number field with ring of integers $\Op_F$. Inspired by the work of Tate \cite{Ta76}, Browkin, Keune and Kolster have made extensive studies in comparing the $K_2$-group $K_2(\Op_F)/p^n$ with the $S$-class group $\Cl_S(F(\mu_{p^n}))$ of $F(\mu_{p^n})$ in a series of papers \cite{Brow92, Keu, Ko89, Ko92}. We should perhaps mention that Coates \cite{C72} has also performed a study of this sort but instead over an (infinite) cyclotomic $\Zp$-extension (also see the related works of Gras \cite{Gras86} and Jaulent \cite{Jau90}).

In this paper, we are interested in extending the works of Browkin, Keune and Kolster to the higher even $K$-groups of $\Op_F$. In slightly more detail, let $i$ be a positive integer. Write $a:=a(F)$ for the largest integer such that $F(\mu_p) = F(\mu_{p^{a}})$, and set $b=b(i)$ to be the largest integer such that $p^b$ divides $i$. Our main theorems consist of comparing the group $K_{2i}(\Op_F)/p^n$ with certain eigenspaces of the $S$-class group of an appropriate cyclotomic extension of $F$. In particular, if $n\leq a+b$, we compare
$K_{2i}(\Op_F)/p^n$ with a certain eigenspace of the $S$-class group of $F(\mu_p)$ (see Theorem \ref{p^n rank thm for small n} for the precise statement). When $n> a+b$, the comparison is done over the field $F(\mu_{p^{n-b}})$ (see Theorem \ref{p^n rank thm}). Note that in the case of $K_2$, the quantity $b$ is always zero and so the comparison for large $n$ is always done over $F(\mu_{p^{n}})$. However, for a higher even $K$-group, the presence of $b$ requires extra care during our discussion.

As illustrations of our results, we specialize them to the context of a quadratic field, where the results become slightly more explicit. In the context of $K_2$, such $p^n$-rank results have also been worked out by Browkin \cite{Brow92} and Zhou \cite{Zhou}. Our discussion here thus generalizes these prior results. We note that many authors have applied the results of Browkin and Keune to obtain explicit results for $K_2$ in various specific number fields (for instance, see \cite{Ch, LZDW, SX, Qin, WZ, Zhou07,Zhou10,ZhouLiang}, where this list is far from being exhaustive). It would definitely be an interesting project to apply the results of this paper to perform analogous studies on the higher even $K$-groups for these said number fields. In this note, we shall only contend ourselves with the case of a quadratic field (although we hope to return to this subject in future works).

We next turn our attention to the case of a totally real field, where we give an algebraic $K$-theoretical formulation of the Kummer's criterion in the sense of Greenberg \cite{Gr73} and Kida \cite{Ki} (see Theorem \ref{kummer} and Remark \ref{kummer remark} for the details). We then prove this via the algebraic $K$-theoretical results developed in this paper. In particular, our proof \textit{does not} make use of the $p$-adic $L$-function.

We end the introductory section giving an outline of the paper. In Section \ref{Alg preliminaries}, we collect certain basic algebraic notion and notation required for our discussion. Section \ref{Arithmetic preliminaries} is where we recall the relation between
the Sylow $p$-subgroups of the even $K$-groups and various Galois cohomology groups. This paves the way for us to prove our main results in Section \ref{Main results section}. Section \ref{quad fields} is where we specialize our main results to the context of a quadratic field. Finally, in Section \ref{kummer section}, we give the algebraic $K$-theoretical formulation of the Kummer's criterion of a totally real field and its proof.

\subsection*{Acknowledgement}
The author would like to thank Georges Gras for his interest and comments on the paper. He is also grateful to the referees for their many valuable comments and suggestions.
  This research is supported by the
National Natural Science Foundation of China under Grant No. 11771164.

\section{Algebraic preliminaries} \label{Alg preliminaries}

As a start, we introduce certain terminology and notation that will be used throughout.
Let $p$ be a prime. For a positive integer $t$, write $v_p(t)$ for the $p$-adic valuation of $t$. In particular, if $t$ is a $p$-power, we have $t = p^{v_p(t)}$. For a finite abelian group $N$, denote by $N[p^n]$ the subgroup of $N$ consisting of elements annihilated by $p^n$. Plainly, one has $N[p^{n-1}]\subseteq N[p^n]$ for $n\geq 1$, where $N[p^0]$ is understood to be the trivial group. The quotient module $N[p^{n}]/N[p^{n-1}]$ can be viewed as a $\mathbb{F}_p$-vector space, where $\mathbb{F}_p$ is the finite field with $p$ elements. Therefore, it makes sense to define the $p^n$-rank $r_{p^n}(N) = \dim_{\mathbb{F}_p}\big(N[p^{n}]/N[p^{n-1}]\big)$.

Throughout the paper, we shall also write $N/p^n$ for $N/p^nN$. The following elementary lemma is left to the reader as an exercise.

\bl \label{alg lemma}
Let $N$ be a finite abelian group. Then we have
\[  r_{p^n}(N) = r_{p^n}(N/p^n).\]
\el

In this paper, we are interested in studying the quantity $r_{p^n}\big(K_{2i}(\Op_F)\big)$ which, by the preceding lemma, leads us to analyse $r_{p^n}\big(K_{2i}(\Op_F)/p^n\big)$.
We shall frequently make use of this viewpoint in the paper without any further mention.

Finally, for a given profinite group $G$ and a $G$-module $M$, we let $M^G$ be the subgroup of $M$ consisting of elements fixed by $G$ and $M_G$ be the largest quotient of $M$ on which $G$ acts trivially. If $M$ is a discrete $G$-module, we write $H^k (G, M)$ for the
$k$-th cohomology group of $G$ with coefficients in $M$.

\section{Arithmetic preliminaries} \label{Arithmetic preliminaries}

For the discussion of our main results, we need to recall the relation between the algebraic $K$-groups and \'etale/Galois cohomology. This relation is essentially well-known (for instances, see \cite{Sou, Vo, WeiKbook}), but for the convenience of the readers, we shall collect some aspects of the said relation in this section. Again, we emphasize that we will always work under the standing assumption that the prime $p$ is odd.

\subsection{$i$-fold tensor product}
To begin with, fix once and for all an algebraic closure $\bar{\Q}$ of $\Q$. Therefore,
an algebraic (possibly infinite) extension of $\Q$ will mean an subfield of our fixed algebraic closure $\bar{\Q}$.
In particular, a number field is understood to be a finite extension of $\Q$ contained in $\bar{\Q}$.

Let $F$ be a number field, whose ring of integers is in turn denoted by $\Op_F$. Throughout the paper, $S$ will always denote the finite set of primes of $F$ which consists precisely of the primes above $p$ and the infinite primes. We then write $\Op_{F,S}$ for the ring of $S$-integers. Let $F_S$ be the maximal algebraic extension of $F$ unramified outside $S$ and denote by $G_S(F)$ the Galois group $\Gal(F_S/F)$.
For a given finite extension $L$ of $F$ contained in $F_S$, we shall write $S(L)$ for the set of primes of $L$ above $S$. Equivalently, $S(L)$ is the set of primes of $L$ which consists of the primes above $p$ and the infinite primes. Since $L$ is contained in $F_S$, it follows from a straightforward verification that $F_S = L_{S(L)}$. We then write $G_S(L)$ for the Galois group $\Gal(F_S/L)$.

Now, a $p$-primary discrete $G_S(F)$-module $M$ may also be viewed as a discrete $G_S(L)$-module via restriction of action. From which, we have restriction map
\[ H^k(G_S(F), M) \lra H^k(G_S(L), M)\]
and corestriction map
\[ H^k(G_S(L), M) \lra H^k(G_S(F), M)\]
on cohomology. For the subsequent discussion of the paper, we require the following standard facts which are recorded here for convenience.

\bl \label{cohomology lemmas}
Suppose that $M$ is a $p$-primary discrete $G_S(F)$-module. Let $L$ be a Galois extension of $F$ contained in $F_S$. Then the following statements are valid.
\begin{enumerate}
    \item[(i)] The corestriction map induces an isomorphism
\[ H^2(G_S(L), M)_{\Gal(L/F)} \cong H^2(G_S(F), M) \]
  \item[(ii)] If the group $\Gal(L/F)$ is finite of order coprime to $p$, then the restriction map induces an isomorphism
   \[  H^2(G_S(F), M) \cong H^2(G_S(L), M)^{\Gal(L/F)}. \]
\end{enumerate}
\el

\bpf
The first assertion follows from reading off the initial term of the Tate spectral sequence (cf. \cite[Theorem 2.5.3]{NSW}). On the other hand, in view of $\Gal(L/F)$ being finite with order coprime to $p$, the Hochschild-Serre spectral sequence
\[ H^r(\Gal(L/F), H^s(G_S(L), M))\Longrightarrow H^{r+s}(G_S(F), M) \]
is concentrated on the line $r=0$ yielding
\[  H^s(G_S(F), M) \cong H^s(G_S(L), M)^{\Gal(L/F)} \]
for every $s\geq 0$. In particular, taking $s=2$, we have the second assertion.
\epf

We now introduce the discrete modules that will be frequently considered in this paper.
Denoting by $\mu_{p^n}$ the cyclic group generated by a primitive $p^n$th-root of unity, we then write $\mu_{p^\infty}$ for the direct limit of the groups $\mu_{p^n}$. These have natural $G_S(F)$-module structures. Furthermore, for an integer $i\geq 2$, the $i$-fold tensor products $\mu_{p^n}^{\otimes i}$ and $\mu_{p^\infty}^{\otimes i}$ can be endowed with $G_S(F)$-module structures via the diagonal action. Therefore, we may speak of the Galois cohomology groups  $H^{k}\big(G_S(F), \mu_{p^n}^{\otimes i}\big)$ and  $H^{k}\big(G_S(F), \mu_{p^\infty}^{\otimes i}\big)$, noting that
\[H^{k}\big(G_S(F), \mu_{p^\infty}^{\otimes i}\big)\cong \ilim_n H^{k}\big(G_S(F), \mu_{p^n}^{\otimes i}\big). \]
In particular, for $1\leq m\leq \infty$, we shall write $\mu_{p^m}^{\otimes i}(L)$ for $H^0(G_S(L),\mu_{p^n}^{\otimes i})= (\mu_{p^n}^{\otimes i})^{G_S(L)}$.

Let $S_p$ be the set of primes of $F$ above $p$. For a finite extension $L$ of $F$ contained in $F_S$, we will write $S_p(L)$ for the set of primes of $L$ above $S_p$. If $w$ is a prime in $S_p(L)$, we then write $L_w$ for the completion of $L$ at $w$. One has a natural group homomorphism
\[ \Gal(\overline{L_w}/L_w) \lra \Gal(\overline{L}/L) \lra \Gal(F_S/L).\]
Via this homomorphism, the $i$-fold tensor product $\mu_{p^m}^{\otimes i}$ may be viewed as a $\Gal(\overline{L_w}/L_w)$-module for $1\leq m\leq \infty$, and we then set $\mu_{p^m}^{\otimes i}(L_w)$  to be $H^0(\Gal(\overline{L_w}/L_w) ,\mu_{p^m}^{\otimes i})$.

The next lemma gives a precise description of these groups.

\bl \label{i-fold invariant}
Suppose that $i\geq 1$. Let $L$ be a finite extension of $F$ contained in $F_S$. Let $\mathcal{L}$ denote either $L$ or $L_w$ for some $w\in S_p(L)$. Set $a(\mathcal{L})$ to be the largest integer such that $\mathcal{L}(\mu_p)$ contains a primitive $p^{a(\mathcal{L})}$th root of unity and set $b(i)=v_p(i)$.
Then we have
\[ \mu_{p^{\infty}}^{\ot i}(\mathcal{L}) = \left\{
                                   \begin{array}{ll}
                                     \mu_{p^{a(\mathcal{L})+b(i)}}^{\ot i}, & \hbox{if $i \equiv 0$ mod $|\mathcal{L}(\mu_p):\mathcal{L}|$,} \\
                                     1, & \hbox{if $i \not\equiv 0$ mod $|\mathcal{L}(\mu_p):\mathcal{L}|$.}
                                   \end{array}
                                 \right.
 \]
In particular, one has
\[ \mu_{p^n}^{\ot i}(\mathcal{L})= \left\{
                   \begin{array}{ll}
                     \mu_{p^n}^{\ot i}, & \hbox{if $i \equiv0$ mod $|\mathcal{L}(\mu_p):\mathcal{L}|$ and $n\leq a(\mathcal{L})+b(i)$,} \\
                     \mu_{p^{a(\mathcal{L})+b(i)}}^{\ot i}, & \hbox{if $i \equiv0$ mod $|\mathcal{L}(\mu_p):\mathcal{L}|$ and $n> a(\mathcal{L})+b(i)$,} \\
               1, & \hbox{if $i \not\equiv 0$ mod $|\mathcal{L}(\mu_p):\mathcal{L}|$.}
                   \end{array}
                 \right.
\]
 \el

\bpf
See \cite[Chap. VI, Proposition 2.2]{WeiKbook}.
\epf

\br
Note that the final equality is also saying that whenever $i =0$ mod $|\mathcal{L}(\mu_p):\mathcal{L}|$ and $n\leq a(\mathcal{L})+b(i)$, $\mu_{p^n}^{\ot i}$ is a trivial $\mathcal{G}$-module, where $\mathcal{G}$ is $G_S(L)$ or $\Gal(\overline{L_w}/L_w)$ accordingly to $\mathcal{L}$ being $L$ or $L_w$.
\er

Let $L$ be a finite extension of $F$. To ease the notation slightly, we shall write $\Cl_{S}(\Op_{L})$ for the $S(L)$-class group $\Cl_{S(L)}(\Op_{L})$. We then write $A_L$ (resp.\ $A_L^S$) for the Sylow $p$-subgroup of $\Cl(\Op_{L})$ (resp., Sylow $p$-subgroup of $\Cl_{S}(\Op_{L})$). The following exact sequence relates this class group with the cohomology group $H^2(G_S(L), \mu_{p^n})$.

\bp \label{Poitou-Tate class group}
Let $L$ be a finite extension of $F$ contained in $F_S$. Then for every $n\geq 1$, one has the following exact sequence
\[ 0\lra A_L^S/p^n \lra H^2(G_S(L), \mu_{p^n})\lra \bigoplus_{w\in S_p(L)}\Z/p^n\Z \lra \Z/p^n\Z \lra 0.\]
\ep

\bpf
See \cite[Satz 4]{Sch} or \cite[Proposition 8.3.11]{NSW}.
\epf

\subsection{A brief interlude on \'etale cohomology}\label{etale Gal}

On the other hand,  we can view $\mu_{p^n}^{\otimes i}$ as an \'etale sheaf over the scheme $\mathrm{Spec}(\Op_{F,S})$ in the sense of \cite[Chap.\ II]{Mi}, and consider the \'etale cohomology groups $H^{k}_{\acute{e}t}\big(\mathrm{Spec}(\Op_{F,S}), \mu_{p^n}^{\otimes i}\big)$.
By \cite[Chap.\ II, Proposition 2.9]{Mi}, the latter is isomorphic to the Galois cohomology groups $H^{k}\big(G_S(F), \mu_{p^n}^{\otimes i}\big)$. Writing $\Zp(i) = \plim_n  \mu_{p^n}^{\otimes i}$ and taking inverse limit, we therefore obtain
\[\plim_n H^{k}_{\acute{e}t}\big(\mathrm{Spec}(\Op_{F,S}), \mu_{p^n}^{\otimes i}\big) \cong \plim_n H^{k}\big(G_S(F), \mu_{p^\infty}^{\otimes i}\big) \cong H^{k}_{\mathrm{cts}}\big(G_S(F), \Zp(i)\big),\]
where $H^{k}_{\mathrm{cts}}(~,~)$ is the continuous cohomology group of Tate (see \cite[Chap 2, \S7]{NSW}), and where the second isomorphism is a consequence of \cite[Corollary 2.7.6 and Theorem 8.3.20(i)]{NSW}.
To simplify notation, we shall write $H^{k}_{\acute{e}t}\big(\Op_{F,S}, \Zp(i)\big)= \plim_n H^{k}_{\acute{e}t}\big(\mathrm{Spec}(\Op_{F,S}), \mu_{p^n}^{\otimes i}\big)$.

\subsection{Algebraic $K$-theory} \label{Alg K-theory}

We come to the algebraic $K$-theoretical aspects. As before, $p$ will always denote an odd prime.
For a ring $R$ with identity, $K_n(R)$ will always denote the algebraic $K$-groups of $R$ in the sense of Quillen \cite{Qui73a, Qui73b} (also see \cite{Kol, WeiKbook}). In particular, we are interested in understanding $K_{2i}(\Op_F)$ for $i\geq 1$, where it is well-known that these groups are finite by the results of Quillen (see \cite{Qui73b, Bo, Gar}).

In \cite{Sou}, Soul\'e connected the higher $K$-groups with \'etale cohomology groups via the $p$-adic Chern class maps
\[ \mathrm{ch}_{i,k}^{(p)}: K_{2i+2-k}(\Op_F)\ot \Zp \lra H^k_{\acute{e}t}\left(\Op_{F,S}, \Zp(i+1)\right)\]
for $i\geq 1$ and $k =1,2$. (For the precise definition of these Chern class maps, we refer readers to loc.\ cit.) The famed Quillen-Lichtenbaum Conjecture asserted that these maps are isomorphisms (for instance, see \cite{Kol, WeiKbook} for the history). Thanks to the gallant efforts of many, this prediction is now known to be true.

\bt \label{chern iso}
For an odd prime $p$, the $p$-adic Chern class maps are isomorphisms for $i\geq 1$ and $k =1,2$.
 \et

\bpf
Soul\'e first proved that these maps are surjective (see \cite[Th\'eor\`{e}me 6(iii)]{Sou}; also see the work of Dwyer and Friedlander \cite[Theorem 8.7]{DF}). It is folklore (for instance, see \cite[Theorem 2.7]{Kol}) that the asserted bijectivity is a consequence of the norm residue isomorphism theorem (previously also known as the Bloch-Kato(-Milnor) conjecture; see \cite{BK, Mil} for details). This said norm residue isomorphism theorem was proved by the deep work of Rost and Voevodsky \cite{Vo} (also see \cite{RW}) with the aid of a patch from Weibel \cite{Wei09}.
\epf

\bc \label{K2 = H2}
Let $p$ be an odd prime. For $i\geq 1$, we have
\[ K_{2i}(\Op_F)[p^\infty] \cong H^2_{\mathrm{cts}}\big(G_{S}(F), \Zp(i+1)\big)\]
Furthermore, for each $n\geq 1$, we have
\[ K_{2i}(\Op_F)/p^n \cong H^2\big(G_{S}(F), \mu_{p^n}^{\ot i+1}\big).\]
\ec

\bpf
Since the group $K_{2i}(\Op_F)$ is finite, it follows that $K_{2i}(\Op_F)[p^\infty]\cong K_{2i}(\Op_F)\ot\Zp$. By the preceding theorem, the latter is isomorphic to $H^2_{\acute{e}t}\big(\Op_{F,S}, \Zp(i+1)\big)$ which identifies with the corresponding continuous Tate cohomology group as noted in Subsection \ref{etale Gal}. The final isomorphism now follows from considering the long exact sequence of the $G_S(F)$-cohomology of the short exact sequence
\[0\lra \Zp(i+1)\stackrel{p^n}{\lra} \Zp(i+1)\lra \mu_{p^n}^{\ot i+1} \lra 0 \]
and noting that $H^3_{\mathrm{cts}}\big(G_{S}(F), -\big) = 0$ (cf. \cite[Proposition 10.11.3]{NSW}). \epf

\br
We should mention that the conclusions of Theorem \ref{chern iso} and Corollary \ref{K2 = H2} are not true in general when $p=2$. As the current paper is concerned with the case of an odd prime $p$, we shall say nothing more on this but refer the interested readers to \cite{RW} for details.
\er

\section{Main results} \label{Main results section}

\subsection{Twist}

As before, $p$ is an odd prime.
Let $m$ be a positive integer. Once and for all, we fix a primitive $p^m$th root of unity $\ze$ in $\bar{\Q}$. Under such a choice, denote by $\sigma_c$ the element of $\Gal(\Q(\mu_{p^m})/\Q)$ which is defined by $\sigma_c(\ze) = \ze^c$. This assignment in turn induces an isomorphism $\Gal(\Q(\mu_{p^m})/\Q)\cong (\Z/p^m\Z)^{\times}$ of groups.

For a number field $F$, we let $\Delta$ denote the maximal subgroup of $\Gal(F(\mu_{p^m})/F)$ with order coprime to $p$. In other words, $\Delta \cong \Gal(F(\mu_p)/F)$. (Of course, the group $\Delta$ may be trivial, for instance, if $F$ contains $\mu_p$.) We shall write $d=|\Delta|$. The Galois group $\Gal(F(\mu_{p^m})/F)$ may be identified with a subgroup of $\Gal(\Q(\mu_{p^m})/\Q)$, and under this identification, the group $\Delta$ can be identified as a subgroup of $(\Z/p\Z)^{\times}$. For our subsequent discussion, we require an explicit description of this embedding of $\Delta$ into $(\Z/p\Z)^{\times}$. Fix a generator $g$ of $(\Z/p^m\Z)^{\times}$, and write $\si: = (\si_g)^{p^{m-1}}$. Under these settings, the elements of the group $\Delta$ may be identified with elements of the form $\si^{k(p-1)/d}$, where $k=0,1,..., d-1$.  Under the surjection $(\Z/p^m\Z)^{\times}\twoheadrightarrow (\Z/p\Z)^{\times}$, the element $g$ is mapped to a generator of $(\Z/p\Z)^{\times}$ which by abuse of notation is still denoted as $g$. Therefore, correspondingly, we have $g^{(p-1)/d}$ being a generator of the subgroup of $(\Z/p\Z)^{\times}$ which is the image of $\Delta$ under the embedding of $\Delta$ into $(\Z/p\Z)^{\times}$.

Let $\om$ be the Teichm\"{u}ller character of the group $(\Z/p\Z)^{\times}$. For a given integer $j$, we define
\[ \e_j:=\e_j(F): = \frac{1}{d}\sum_{k=0}^{d-1}\om(g)^{jk(p-1)/d}\sigma^{-k(p-1)/d}\]
which lives in $\Zp[\Delta]\subseteq \Zp[\Gal(F(\mu_{p^m})/F)]$. One can check easily that $\e_j=\e_{j'}$ whenever $j \equiv j'$ (mod $|\Delta|$), and that $\e_0,...,\e_{d-1}$ forms a collection of primitive idempotents of  the group ring $\Zp[\Delta]$.

\bl \label{mu twist}
Let $M$ be a $\Zp[\Delta]$-module. Then the following statements are valid.
\begin{enumerate}
  \item[$(i)$] $M^\Delta = \e_0M$.
  \item[$(ii)$] For every $i\geq 1$, we have  $\big(\mu_{p^m}^{\ot i}\ot M\big)^\Delta = \mu_{p^m}^{\ot i}\ot \e_{-i}M$.
\end{enumerate}
\el

\bpf
The first assertion is immediate from the definition of $\e_0$. For the second, it suffices to show that
\[ \e_{j+1}(\mu_{p^m}\ot M) = \mu_{p^m}\ot \e_j M\]
for every $j$. Indeed, the conclusion of the lemma follows from a recursive application of the said equality. Via the natural surjection $(\Z/p^m\Z)^{\times}\twoheadrightarrow (\Z/p\Z)^{\times}$, we may also view $\om$ as a character of $(\Z/p^m\Z)^{\times}$. It then follows that $\om(g)^{p-1}=1$ and $\om(g)=g$ mod $p$, where the latter in turn implies that $\om(g)=\om(g)^{p^{m-1}} = g^{p^{m-1}}$ mod $p^{m-1}$. Consequently, we have
\[ \si(\ze) = \ze^{g^{p^{m-1}}} = \ze^{\om(g)}. \]
Let $x\in M$. Then we have
\[ \e_{j+1}(\ze\ot x) = \frac{1}{d}\sum_{k=0}^{d-1}\om(g)^{(j+1)k(p-1)/d}\sigma^{-k(p-1)/d}(\ze\ot x)\]
 \[= \frac{1}{d}\sum_{k=0}^{d-1}\om(g)^{jk(p-1)/d}\sigma^{-k(p-1)/d}(\ze^{\om(g)^{k(p-1)/d}}\ot x) \]
\[= \frac{1}{d}\sum_{k=0}^{d-1}\om(g)^{jk(p-1)/d}(\ze\ot \sigma^{-k(p-1)/d}x) = \ze\ot \e_jx. \]
This establishes our claim and the proof of the lemma is thus complete.
\epf

In what follows, we will study the $p^n$-rank of $K_{2i}(\Op_F)$ for a given number field $F$ and integer $i\geq 1$. Let $a:=a(F)$ be the largest integer such that $F(\mu_p) = F(\mu_{p^{a}})$. Set $b=v_p(i)$. In other words, $b$ is the largest integer such that $p^b$ divides $i$. We first consider the case when $n\leq a+b$.

\subsection{$p^n$-rank for $n\leq a+b$}
 For a given integer $i\geq 1$, let $S_p^{(i)}$ be the set of primes $v$ in $S_p$ such that $i$ is divisible by $|\Delta_v|$, where $\Delta_v$ is the decomposition group of $\Delta$ at $v$. In particular, $S_p^{(1)}$ is the set of primes $v$ in $S_p$ which splits completely in $F(\mu_p)/F$. We can now state and prove the following.

\bt \label{p^n rank thm for small n} Suppose that $n\leq a+b$.
\begin{itemize}
  \item[$(i)$] If $i\equiv 0$ mod $|F(\mu_p):F|$, we have the following exact sequence
\[ 0\lra \mu_{p^n}^{\ot i}\ot A^S_F\lra  K_{2i}(\Op_F)/p^n \lra \bigoplus_{v\in S_p}\mu_{p^n}^{\ot i} \lra \mu_{p^n}^{\ot i} \lra 0.\]
  \item[$(ii)$] If $i\not\equiv 0$ mod $|F(\mu_p):F|$, we have the following exact sequence
\[ 0\lra \mu_{p^n}^{\ot i}\ot \e_{-i}A^S_{F(\mu_p)} \lra  K_{2i}(\Op_F)/p^n \lra \bigoplus_{v\in S_p^{(i)}}\mu_{p^n}^{\ot i} \lra 0.\]
\end{itemize}
\et

\bpf
We first consider the situation when  $i\equiv 0$ mod $|F(\mu_p):F|$.  From Proposition \ref{Poitou-Tate class group}, we have the following exact sequence
\[ 0\lra A^S_F/p^n \lra H^2(G_S(F), \mu_{p^n})\lra \bigoplus_{v\in S_p}\Z/p^n\Z \lra \Z/p^n\Z \lra 0.\]
Applying $\mu_{p^n}^{\ot i}\ot-$, we obtain
\[ 0\lra \mu_{p^n}^{\ot i}\ot A^S_F\lra \mu_{p^n}^{\ot i}\ot H^2(G_S(F), \mu_{p^n})\lra \bigoplus_{v\in S_p}\mu_{p^n}^{\ot i} \lra \mu_{p^n}^{\ot i} \lra 0.\]
Since $i\equiv 0$ mod $|F(\mu_p):F|$, the group $G_S(F)$ acts trivially on $\mu_{p^n}^{\ot i}$, and so we have
\[  \mu_{p^n}^{\ot i}\ot H^2(G_S(F), \mu_{p^n}) \cong  H^2(G_S(F), \mu_{p^n}^{\ot i+1}).\]
By Corollary \ref{K2 = H2}, the latter is $K_{2i}(\Op_F)/p^n$. This proves (i).

Now suppose that $i\not\equiv 0$ mod $|F(\mu_p):F|$. In this case, we have $G_S(F(\mu_p))$ acting trivially on $\mu_{p^n}^{\ot i}$ by Lemma \ref{cohomology lemmas}, since $a(F(\mu_p))=a$. By a similar argument as above, we obtain the following exact sequence
\[ 0\lra \mu_{p^n}^{\ot i}\ot A^S_{F(\mu_p)}\lra  H^2(G_S(F(\mu_p)), \mu_{p^n}^{\ot i+1})\lra \bigoplus_{w\in S_p(F(\mu_p))}\mu_{p^n}^{\ot i} \lra \mu_{p^n}^{\ot i} \lra 0.\]
Taking $\Delta$-invariant, we have
\begin{multline*}
  0\lra \big(\mu_{p^n}^{\ot i}\ot A^S_{F(\mu_p)})\big)^\Delta\lra  H^2(G_S(F(\mu_p)), \mu_{p^n}^{\ot i+1})^\Delta  \\
\lra \Big(\bigoplus_{w\in S_p(F(\mu_p))}\mu_{p^n}^{\ot i}\Big)^\Delta \lra \big(\mu_{p^n}^{\ot i}\big)^\Delta \lra 0.
\end{multline*}
Note that the above sequence is still exact, as $|\Delta|$ is coprime to $p$. In view of this latter coprime property, one may also apply Lemma \ref{cohomology lemmas} to conclude that
\[ H^2\big(G_S(F(\mu_p)), \mu_{p^n}^{\ot i+1}\big)^\Delta\cong H^2(G_S(F), \mu_{p^n}^{\ot i+1}),\]
where the latter is $K_{2i}(\Op_F)/p^n$ by Corollary \ref{K2 = H2}. The first term in the exact sequence is \[\mu_{p^n}^{\ot i}\ot \e_{-i}A^S_{F(\mu_p)}\] by Lemma \ref{mu twist}.
On the other hand, as $i\not\equiv 0$ mod $|F(\mu_p):F|$, one has $\big(\mu_{p^n}^{\ot i}\big)^\Delta =0$ by Lemma \ref{cohomology lemmas}. For the remaining local terms, we note that
\[ \Big(\bigoplus_{w\in S_p(F(\mu_p))}\mu_{p^n}^{\ot i}\Big)^\Delta = \bigoplus_{v\in S_p}\Big(\bigoplus_{w \mid v}\mu_{p^n}^{\ot i}\Big)^\Delta \cong \bigoplus_{v\in S_p}\big(\mu_{p^n}^{\ot i}\big)^{\Delta_v}, \]
where the second isomorphism follows from Shapiro's lemma.
By another application of Lemma \ref{cohomology lemmas}, one then has
\[(\mu_{p^n}^{\ot i}\big)^{\Delta_v} = \left\{
                                          \begin{array}{ll}
                                            \mu_{p^n}^{\ot i}, & \hbox{if $i$ is divisible by $|\Delta_v|$,} \\
                                            0, & \hbox{otherwise.}
                                          \end{array}
                                        \right.
 \]
Combining all the above observations, we obtain the conclusion of (ii).
\epf

We discuss some consequences of the preceding theorem.

\bc \label{p-rank cor small}
Suppose that $n\leq a+b$. Then \[ r_{p^n}\big( K_{2i}(\Op_F)\big) =\left\{
                                      \begin{array}{ll}
                                        r_{p^n}\big( A^S_F \big) +|S_p|-1, & \hbox{if $i\equiv 0$ mod $|F(\mu_p):F|$,} \\
                                        r_{p^n}\big( \e_{-i}A^S_{F(\mu_p)} \big) +|S_p^{(i)}|, & \hbox{if $i\not\equiv 0$ mod $|F(\mu_p):F|$.}
                                      \end{array}
                                    \right.
 \]
\ec

\bpf Suppose that $i=0$ mod $|F(\mu_p):F|$. Set $N$ to be the kernel of the map $\bigoplus_{v\in S_p}\mu_{p^n}^{\ot i} \lra \mu_{p^n}^{\ot i}$. Then it follows from Theorem \ref{p^n rank thm for small n}(i) that we have a short exact sequence
\[ 0\lra \mu_{p^n}^{\ot i}\ot A^S_F\lra  K_{2i}(\Op_F)/p^n \lra N \lra 0\]
of $\Z/p^n\Z$-modules. Plainly, $N$ is a free $\Z/p^n\Z$-module of rank $|S_p|-1$, and so the above short exact sequence splits, which in turn supplies an isomorphism
\[K_{2i}(\Op_F)/p^n \cong \big(\mu_{p^n}^{\ot i}\ot A^S_F\big)\oplus N  \]
of $\Z/p^n\Z$-modules. The asserted $p^n$-rank of $K_{2i}(\Op_F)$ is now an immediate consequence of this.

Suppose that $i\neq 0$ $($mod $|\Delta|)$. By Theorem \ref{p^n rank thm for small n}, we have a short exact sequence
 \[ 0\lra \mu_{p^n}^{\ot i}\ot \e_{-i}A^S_{F(\mu_p)} \lra  K_{2i}(\Op_F)/p^n \lra \bigoplus_{v\in S_p^{(i)}}\mu_{p^n}^{\ot i} \lra 0,\]
where $\bigoplus_{v\in S_p^{(i)}}\mu_{p^n}^{\ot i}$ is a free $\Z/p^n\Z$-module of rank $|S_p^{(i)}|$. Therefore, the sequence splits and the asserted $p^n$-rank of $K_{2i}(\Op_F)$ follows from this.
\epf

We record the following special case of the preceding corollary which generalizes a previous result of Kolster (see \cite[Corollary 1.8]{Ko89} and \cite[Corollary 3.4]{Ko92}), where a formula for the $p^n$-rank of $K_2$ is obtained.

\bc
Suppose $n\leq a+b$, $|F(\mu_p):F|=2$ and $i$ is an odd integer. Then we have
\[ r_{p^n} (K_{2i}(\Op_F)) = r_{p^n}\big(\e_{1}A^S_{F(\mu_{p})}\big) + |S_p^{(1)}|.\]
\ec

\subsection{$p^n$-rank for $n>a+b$}

We now consider the situation when $n$ is greater than $a+b$.

\bt \label{p^n rank thm}
If $n> a+b$, we have the following exact sequence
\[ 0\lra \big(\mu_{p^n}^{\ot i}\ot A^S_{E}\big)_{G}\lra  K_{2i}(\Op_F)/p^n \lra \bigoplus_{v\in S_p}\mu_{p^n}^{\ot i}(F_v) \lra \mu_{p^n}^{\ot i}(F) \lra 0\]
where $E:=E_n:=F(\mu_{p^{n-b}})$ and $G=\Gal(E/F)$.
\et

In preparation of the proof of the theorem, we shall first establish the following lemma which generalizes an observation made by Keune in \cite[Lemma 6.5]{Keu}. We note that our approach of proof differs from that of Keune.

\bl \label{coh trivial}
Retain the settings of Theorem \ref{p^n rank thm}. Then $\mu_{p^{n}}^{\ot i}$ is a cohomologically trivial $H$-module for every subgroup $H$ of $G$.
\el

\bpf
By replacing $F$ by $E^H$, it suffices to prove that $\mu_{p^{n}}^{\ot i}$ is a cohomologically trivial $G$-module.
Set $\Ga= \Gal(F(\mu_{p^\infty})/F)$ and $\Ga_n = \Gal(F(\mu_{p^\infty})/E)$. Consider the following spectral sequence
\[  H^r\big(G, H^s(\Ga_n, \mu_{p^\infty}^{\ot i})\big)\Longrightarrow H^{r+s}(\Ga, \mu_{p^\infty}^{\ot i}). \]
Since the groups $\Ga_n$ and $\Ga$ have $p$-cohomological dimension one, one has
\[H^t(\Ga_n, \mu_{p^\infty}^{\ot i})=H^t(\Ga, \mu_{p^\infty}^{\ot i})=0\]
 for $t\geq 2$. On the other hand, Tate Lemma \cite{Ta} tells us that
\[H^1(\Ga_n, \mu_{p^\infty}^{\ot i})=H^1(\Ga, \mu_{p^\infty}^{\ot i})=0.\]
 Hence the spectral sequence degenerates yielding
\[ H^r\big(G, H^0(\Ga_n, \mu_{p^\infty}^{\ot i})\big) =0\]
for every $r\geq 1$. But one also has $ H^0(\Ga_n, \mu_{p^\infty}^{\ot i}) = \mu_{p^{n}}^{\ot i}$ by Lemma \ref{cohomology lemmas}. This completes the proof of the lemma. \epf

We can now give the proof of Theorem \ref{p^n rank thm}.

\bpf[Proof of Theorem \ref{p^n rank thm}]
 In view that $G_S(E)$ acts trivially on $\mu_{p^n}^{\ot i}$, we may apply a similar argument as in the beginning of the proof of Theorem \ref{p^n rank thm for small n} to obtain the following exact sequence
\[ 0\lra \mu_{p^n}^{\ot i}\ot A^S_E\lra  H^2(G_S(E), \mu_{p^n}^{\ot i+1})\lra \bigoplus_{w\in S_p(E)}\mu_{p^n}^{\ot i} \lra \mu_{p^n}^{\ot i} \lra 0.\]
 By Lemma \ref{coh trivial}, we see that $\mu_{p^n}^{\ot i}$ is a cohomologically trivial $G$-module. We now proceed to show that $\bigoplus_{w\in S_p(E)}\mu_{p^n}^{\ot i}$ is also a cohomologically trivial $G$-module. Observe that
\[ H^r\Big(G, \bigoplus_{w\in S_p(E)}\mu_{p^n}^{\ot i}\Big) \cong \bigoplus_{v\in S_p}H^r\Big(G, \bigoplus_{w|v}\mu_{p^n}^{\ot i}\Big) \cong  \bigoplus_{v\in S_p}H^r\Big(G_v, \mu_{p^n}^{\ot i}\Big),\]
where the second isomorphism is a consequence of Shapiro's lemma. Since $G_v$ may be identified as a subgroup of $G$, we may apply Lemma \ref{coh trivial} to conclude that $H^r\Big(G_v, \mu_{p^n}^{\ot i}\Big) =0$ for $r\geq 1$. Hence this shows that $\bigoplus_{w\in S_p(E)}\mu_{p^n}^{\ot i}$ is a cohomologically trivial $G$-module. Set $Y:= \ker \Big(\displaystyle\bigoplus_{w\in S_p(E)}\mu_{p^n}^{\ot i} \lra \mu_{p^n}^{\ot i} \Big)$. A straightforward cohomological consideration shows that $Y$ is also a cohomologically trivial $G$-module.
Therefore, after taking $-_G$ functor, we have the following short exact sequences
\[0\lra \Big(\mu_{p^n}^{\ot i}\ot A^S_E\Big)_G\lra  \Big(H^2(G_S(E), \mu_{p^n}^{\ot i+1})\Big)_G \lra Y_G\lra 0, \]
\[ 0\lra Y_G\lra \Big(\bigoplus_{w\in S_p(E)}\mu_{p^n}^{\ot i}\Big)_G \lra \Big(\mu_{p^n}^{\ot i}\Big)_G \lra 0.\]
Splicing these exact sequences, we obtain
\[0\lra \Big(\mu_{p^n}^{\ot i}\ot A^S_E\Big)_G\lra  \Big(H^2(G_S(E), \mu_{p^n}^{\ot i+1})\Big)_G \lra \Big(\bigoplus_{w\in S_p(E)}\mu_{p^n}^{\ot i}\Big)_G \lra \Big(\mu_{p^n}^{\ot i}\Big)_G \lra 0. \]
By Lemma \ref{cohomology lemmas}, $\Big(H^2(G_S(E), \mu_{p^n}^{\ot i+1})\Big)_G \cong H^2(G_S(F), \mu_{p^n}^{\ot i+1})$ which is precisely $K_{2i}(\Op_F)/p^n$ by Corollary \ref{K2 = H2}. Finally, for every cohomologically trivial $G$-module $N$, one has a natural isomorphism $N_G\cong N^G$ induced by the norm map. Therefore, for the rightmost two terms, we may switch  $-_G$ to $-^G$.
 Putting all these observations together, we obtain the conclusion of our theorem.
\epf

\section{Quadratic fields} \label{quad fields}

We specialize to the quadratic fields. Let $m$ be a squarefree integer $\neq \pm 1, (-1)^{(p-1)/2}p$. Set $F=\Q(\sqrt{m})$. Note that in this sitation, one has $a = a(F)= 1$. As before, we continue to write $b=v_p(i)$.

We begin with the case of small $n$.

\bp
Suppose that $n\leq 1+b$. Then the following statements are valid.
\begin{itemize}
  \item[$(i)$] If $i\equiv 0$ mod $p-1$, then
\[ r_{p^n}\big( K_{2i}(\Op_F)\big) =\left\{
                                      \begin{array}{ll}
                                        r_{p^n}\big( A^S_F \big) +1, & \hbox{if $p$ splits in $F/\Q$,} \\
                                        r_{p^n}\big( A^S_F \big) , & \hbox{otherwise.}
                                      \end{array}
                                    \right.
 .\]
  \item[$(ii)$] If $i\not\equiv 0$ mod $p-1$, then we have
\begin{multline*}
r_{p^n} (K_{2i}(\Op_F)) \\
     =\left\{                              \begin{array}{ll}
                                r_{p^n}\big(\e_{-i}A^S_{F(\mu_p)}\big) +1, & \hbox{if $i\equiv 0$ mod $(p-1)/2$ and $m=(-1)^{(p-1)/2}pm_1$}\\
&\hbox{ with $p$ being split in $\Q(\sqrt{m_1})/\Q$,} \\
                                r_{p^n}\big(\e_{-i}A^S_{F(\mu_p)}\big) , & \hbox{otherwise.}
                              \end{array}
                            \right.
   \end{multline*}

\end{itemize}

\ep

\bpf (i) Since $|S_p| = 2$ or $1$ accordingly to the prime $p$ splitting in $F/\Q$ or not, the asserted $p^n$-rank formula is a consequence of Corollary \ref{p-rank cor small}.

(ii) From Corollary \ref{p-rank cor small}, one has
 \[r_{p^n} (K_{2i}(\Op_F)) = r_{p^n}\big(\e_{-i}A^S_{F(\mu_p)}\big) + |S^{(i)}|.\]
From the discussion in \cite[Section 3]{Zhou}, we see that for a prime $v\in S_p$, $|F_v(\mu_p):F_v| =  (p-1)/2$
precisely when the prime $p$ is both ramified in $F/\Q$ and split in $\Q(\sqrt{m_1})/\Q$, where $m=(-1)^{(p-1)/2}pm_1$. The discussion in \cite[Section 3]{Zhou} also tells us that in all other cases, one has $|F_v(\mu_p):F_v| =  p-1$. Hence we have $|S_p^{(i)}|=1$ when $i$ is divisible by $(p-1)/2$, and $m=(-1)^{(p-1)/2}pm_1$ with $p$ splitting in $\Q(\sqrt{m_1})/\Q$, and $|S_p^{(i)}|=0$ otherwise. This yields the asserted conclusion. \epf

For large $n$, we may apply Theorem \ref{p^n rank thm} to obtain the following. Since the argument is quite similar to that of the preceding proof, we shall omit the proof.

\bp
Suppose that $n> 1+b$. Set $E=F(\mu_{p^{n-b}})$ and $G=\Gal(E/F)$. Then the following statements are valid.
\begin{itemize}
  \item[$(i)$] Suppose that $i\equiv 0$ mod $p-1$.
  \begin{itemize}
    \item[$(a)$] If $p$ is non-split in $F/\Q$, then we have the following exact sequence
\[ \big(\mu_{p^n}^{\ot i}\ot A^S_{E}\big)_{G}\cong  K_{2i}(\Op_F)/p^n \]
and therefore
\[ r_{p^n}\big( K_{2i}(\Op_F)\big) = r_{p^n}\Big(\big(\mu_{p^n}^{\ot i}\ot A^S_E\big)_{G} \Big).\]
\item[$(b)$] If $p$ splits in $F/\Q$, then we have \[ 0\lra \big(\mu_{p^n}^{\ot i}\ot A^S_E\big)_{G}\lra  K_{2i}(\Op_F)/p^n \lra \mu_{p^{1+b}}^{\ot i} \lra 0.\]
\end{itemize}
  \item[$(ii)$] Suppose that $i\not\equiv 0$ mod $p-1$, then we have
\[  K_{2i}(\Op_F)/p^n \cong \big(\mu_{p^n}^{\ot i}\ot A^S_E\big)_{G}\]
except when one has  $i\equiv 0$ mod $(p-1)/2$ and $m=(-1)^{(p-1)/2}pm_1$ with $p$ being split in $\Q(\sqrt{m_1})/\Q$. In this exceptional case, one instead has a short exact sequence
\[ 0\lra \big(\mu_{p^n}^{\ot i}\ot A^S_E\big)_{G}\lra  K_{2i}(\Op_F)/p^n \lra \mu_{p^{1+b}}^{\ot i} \lra 0.\]

\end{itemize}
\ep

\section{On Kummer's criterion} \label{kummer section}

In this final section, $F$ will always denote a totally real field. We let $\ze_p$ denote a fixed primitive $p$-th root of unity, and set $d=|F(\ze_p):F|$. The following theorem is the main result of this section.

\bt \label{kummer}
Suppose that $F$ is a totally real field such that all the primes in $S_p\big(F(\ze_p+\ze_p^{-1})\big)$ do not split in $F(\ze_p)$. Then the following statements are equivalent.
\begin{itemize}
\item[$(1)$] The class number of $F(\ze_p)$ is divisible by $p$.
\item[$(2)$] The prime $p$ divides the order of $K_{2i}(\Op_F)$ for some $1\leq i \leq d-1$ with $i$ being odd.
\end{itemize}
\et

\br \label{kummer remark} Before giving a proof of the above theorem, we explain how the said theorem is equivalent to Kummer's criterion in the sense of \cite[Theorem 1]{Gr73} and \cite[Theorem 1]{Ki}, and thus, can be thought as an algebraic $K$-theoretical formulation of it. Write $|~|_p$ for the $p$-adic norm which is normalized with $|p|_p = 1/p$. Denote by $\ze_{F}(s)$ the Dedekind zeta function of $F$. For $a,b\in \Q-\{0\}$, we write $a\sim_p b$ if  $a/b$ is a $p$-unit. For each odd integer $i$ with $1\leq i \leq d-1$, it follows from the main conjecture of Iwasawa as proven by Wiles \cite{Wiles} that
\[ \ze_{F}\big(1-(i+1)\big) \sim_p \frac{|K_{2i}(\Op_F)|}{|K_{2i+1}(\Op_F)|}. \]
(Strictly speaking, the theorem of Wiles will only yield
\[ \ze_{F}\big(1-(i+1)\big) \sim_p \frac{|H^2(G_S(F), \Zp(i+1)|}{|H^1(G_S(F), \Zp(i+1)|}; \]
see the paper of B\'ayer and Neukirch \cite{BN}. But we now know that this cohomological version is equivalent to the $K$-theoretical version as stated above in view of the work of Rost-Voevodsky \cite{Vo}).

From \cite[Chap.\ VI, Theorem 9.5]{WeiKbook}, we see that $\big|K_{2i+1}(\Op_F)\big|\sim_p w^{(p)}_{i+1}(F)$, where $w_j^{(p)}(F)$ is the order of $\mu_{p^\infty}^{\otimes j}(F)=(\mu_{p^\infty}^{\otimes j})^{\Gal(\bar{F}/F)}$. Lemma \ref{i-fold invariant} then tells us that
\[w^{(p)}_{i+1}(F) = \left\{
            \begin{array}{ll}
              1, & \hbox{$1\leq i <d-1$, $i$ odd} \\
              p^{a(F)}, & \hbox{$i= d-1$.}
            \end{array} \right. \]
Therefore, upon combining these observations, we obtain
\[|K_{2i}(\Op_F)| \sim_p \left\{
            \begin{array}{ll}
              \ze_{F}\big(1-(i+1)\big), & \hbox{$1\leq i <d-1$, $i$ odd} \\
              p^{a(F)}\ze_{F}\big(1-d\big), & \hbox{$i= d-1$.}
            \end{array} \right. \]
Therefore, statement (2) of Theorem \ref{kummer} is equivalent to saying that $p$ divides one of the numerator of the following rational numbers
 \[ \ze_{F}\big(1-(i+1)\big) \quad (\mbox{$1\leq i <d-1$, $i$ odd}), \quad p^{a(F)}\ze_{F}\big(1-d\big). \]
It thus follows that our theorem is equivalent to \cite[Theorem 1]{Gr73} and \cite[Theorem 1]{Ki}.
\er

We proceed with the proof of Theorem \ref{kummer}, where we emphasize that our proof does not make use of $p$-adic $L$-functions. We however should make the following remark.

\br
Although the proof of Theorem \ref{kummer} does not make use of $p$-adic $L$-function, in order to see that the said theorem is equivalent to Kummer's criterion, one requires to be able to relate the special values of the Dedekind zeta function and the size of the $K$-groups, and such relation is predicted by a conjecture of Lichtenbaum \cite{Lic72}. At our current knowledge, it would seem that the only way to study this relation is via the Iwasawa main conjecture, whose formulation itself will require the $p$-adic $L$-functions.
\er

 For the proof of Theorem \ref{kummer}, we require the following lemma.

\bl \label{kummer lemma}
Retain the setting of Theorem \ref{kummer}. For every odd integer $i$ such that $1\leq i\leq d-1$, we have
\[ r_p\big(K_{2i}(\Op_F)\big) = r_p\big(\e_{-i} A_{F(\ze_p)}\big).\]
\el

\bpf
In the proof, we shall write $E= F(\ze_p)$. Recall that by Corollary \ref{p-rank cor small}, one has
\[ r_p\big(K_{2i}(\Op_F)\big) = r_p\big(\e_{-i} A^S_E \big) + |S_p^{(i)}|.\]
From the assumption that the primes in $S_p\big(F(\ze_p+\ze_p^{-1})\big)$ do not split in $E$, we see that $|\Delta_v|$ is even for every $v\in S_p$. But since $i$ is odd, this in turn implies that the set  $S_p^{(i)}$ is empty. It therefore remains to show that
\[ \e_{-i}A^S_E= \e_{-i}A_E\]
for every odd $i$. Since $d$ is even, this is equivalent to showing that $\e_{i}A^S_E= \e_{i}A_E$ for every odd $i$. Let $\la$ be the natural surjection $\Cl(\Op_E)\twoheadrightarrow \Cl_S(\Op_E)$. Then $\ker \la$ is generated by the class of primes in $S_p(E)$. Now observe that for an odd integer $i$, we have

\begin{align*}
 \e_{i} = & ~\frac{1}{d}\sum_{k=0}^{d-1}\om(g)^{-ik(p-1)/d}\sigma^{-k(p-1)/d} \\
   = & ~\frac{1}{d}\sum_{k=0}^{\frac{d}{2}-1}\left(\om(g)^{ik(p-1)/d}\sigma^{-k(p-1)/d} + \om(g)^{i(k+\frac{d}{2})(p-1)/d}\sigma^{-(k+\frac{d}{2})(p-1)/d}\right) \\
 =& ~\frac{1}{d}\sum_{k=0}^{\frac{d}{2}-1}\om(g)^{ik(p-1)/d}\sigma^{-k(p-1)/d}\left( 1 + (-1)^i\sigma^{-\frac{d}{2}(p-1)/d}\right) \\
=& ~\frac{1}{d}\sum_{k=0}^{\frac{d}{2}-1}\om(g)^{ik(p-1)/d}\sigma^{-k(p-1)/d}\left( 1 -\sigma^{\frac{p-1}{2}}\right),
\end{align*}
where we note that $\sigma^{-\frac{p-1}{2}} = \sigma^{\frac{p-1}{2}}$ is the complex conjugation and hence a generator of the group $\Gal(E/F(\ze_p+\ze_p^{-1}))$.
Since the primes in $S_p\big(F(\ze_p+\ze_p^{-1})\big)$ do not split in $F(\ze_p)$, they are invariant under the complex conjugation $\sigma^{\frac{p-1}{2}}$ and so $\ker \la$ is annihilated by $1 -\sigma^{\frac{p-1}{2}}$. This in turn implies that $(\ker \la)[p^\infty]$ is annihilated by $\e_i$ for every odd $i$, which in turn yields $\e_{i}A^S_E= \e_{i}A_E$ for every odd $i$.  The proof of the lemma is now complete.          \epf

\bpf[Proof of Theorem \ref{kummer}]
Let $A$ denote the Sylow $p$-subgroup of the class group of $E:=F(\ze_p)$. Then there is a decomposition
\[ A = \bigoplus_{i=0}^{d-1}\e_i A,\]
where one notes that $\displaystyle\bigoplus_{\substack{0\leq i \leq d-1\\ i~\mathrm{even}}}\e_i A$ is the Sylow $p$-subgroup of the class group of $F(\ze_p+\ze_p^{-1})$.
By \cite[Section 4]{Gr73}, the class number of $E$ is divisible by $p$ if and only the relative class number of $E/F(\ze_p+\ze_p^{-1})$ is divisible by $p$. In view of the above decomposition, this latter divisibility is therefore equivalent to
\[ r_p\left(\bigoplus_{\substack{1\leq i \leq d-1\\ i~\mathrm{odd}}}\e_i A \right)\geq 1.\]
By Lemma \ref{kummer lemma}, this in turn is equivalent to $r_p\big(K_{2i}(\Op_F)\big)\geq 1$ for some odd $i$ with $1\leq i \leq d-1$. The proof of the theorem is therefore complete.
\epf

\end{document}